\documentstyle[12pt,amsfonts,leqno,amscd]{amsart}
\newtheorem{Theorem}{Theorem}[section]
\newtheorem{Definition}[Theorem]{Definition}
\newtheorem{Proposition}[Theorem]{Proposition}

\newtheorem{Lemma}[Theorem]{Lemma}
\newtheorem{Corollary}[Theorem]{Corollary}
\theoremstyle{remark}
\newtheorem{Example}[Theorem]{Example}

\def\il{\int\limits_}

\def\al{\alpha}
\def\gm{\gamma}
\def\th{\theta}

\def\Gm{\Gamma}
\def\Th{\Theta}

\def\sm{\setminus}

\def\re{{\mathbf {Re\,}}}
\def\im{{\mathbf {Im\,}}}
\def\be{\begin{enumerate}}
\def\ee{\end{enumerate}}
\def\bT{\begin{Theorem}}
\def\eT{\end{Theorem}}
\def\bP{\begin{Proposition}}
\def\eP{\end{Proposition}}
\def\bD{\begin{Definition}}
\def\eD{\end{Definition}}
\def\bE{\begin{Example}}
\def\eE{\end{Example}}
\def\bL{\begin{Lemma}}
\def\eL{\end{Lemma}}
\def\bC{\begin{Corollary}}
\def\eC{\end{Corollary}}

\newcommand{\ra}{\rightarrow}
\begin{document}
\title{Smooth submanifolds intersecting any analytic curve in
a discrete set}
\author{Dan Coman, Norman Levenberg and Evgeny A.
Poletsky}\date{November 2003} 
\thanks{D. Coman and E. A. Poletsky were supported by NSF grants}
\address{D. Coman and E. A. Poletsky: dcoman@@syr.edu, eapolets@@syr.edu,
Department of Mathematics,  215 Carnegie Hall, Syracuse
University, Syracuse, NY 13244, USA}
\address{N. Levenberg: levenber@@math.auckland.ac.nz, Department of
Mathematics, University of Auckland, Private Bag 92019 Auckland,
New Zealand}
\begin{abstract} We construct examples of $C^\infty$ smooth
submanifolds in ${\Bbb C}^n$ and ${\Bbb R}^n$ of codimension 2 and
1, which intersect every complex, respectively real, analytic
curve in a discrete set. The examples are realized either as
compact tori or as properly imbedded Euclidean spaces, and are the
graphs of quasianalytic functions. In the complex case, these
submanifolds contain real $n$-dimensional tori or Euclidean spaces
that are not pluripolar while the intersection with any complex
analytic disk is polar.
\end{abstract}
\maketitle
\section{Introduction}
If a real analytic $(2n-2)$-dimensional submanifold $R$ in ${\Bbb
C}^n$ either intersects every complex analytic disk in a discrete
set or contains the disk, then $R$ is a complex submanifold. A
natural question arises: is this true when $R$ is merely smooth?
\par The main result of this paper is:
\bT\label{T:tor} There exists a smooth, compact manifold $R$ in ${\Bbb
C}^n$, diffeomorphic to a $(2n-2)$-dimensional torus, which
intersects every analytic disk in a discrete set. Moreover, $R$
contains a smooth submanifold $M$, diffeomorphic to an
$n$-dimensional torus, which is not pluripolar.\eT
\par Since $R$ does not contain any analytic disk, it is not a
complex manifold. Thus in the category of smooth manifolds the
discreteness of intersections with analytic disks does not imply
that a submanifold is complex.
\par To explain the last statement of the theorem,
we recall that a set $K$ in ${\Bbb C}^n$ is {\it pluripolar} if
there is a plurisubharmonic function $u\not \equiv -\infty$ which
is equal to $-\infty$ on $K$. If $n=1$, pluripolar sets are
classical polar sets. In general, it is hard to detect whether a
set is pluripolar.
\par In 1971 L. I. Ronkin introduced in \cite{R} the notion of
$\Gm$-capacity of a set, which is computed in terms of the
capacities of intersections of that set with complex lines. Ch.
Kiselman \cite{Ki} and, independently, A. Sadullaev \cite{Sa}
constructed a real algebraic surface in ${\Bbb C}^2$ which
intersects any complex line in at most 4 points. So the
$\Gm$-capacity of this set is 0, while the set is not pluripolar.
\par The question whether the polarity of intersections of a set
implies the pluripolarity of that set when complex lines are
replaced by one-dimensional varieties of higher degree was open
since that time. It was posed as an open problem by E. Bedford in
his survey \cite{B} as follows:
\par {\sl Let $E\subset {\Bbb C}^n$. Suppose $E\cap V$ is polar
in $V$ for each germ $V$ of an irreducible, one-dimensional
complex variety $V\subset {\Bbb C}^n$. Is $E$ pluripolar in ${\Bbb
C}^n$}?
\par So the last part of our theorem answers this question in the
negative.
\par The submanifolds in our examples are graphs of quasianalytic
functions (see Section 2) that belong to Denjoy--Carleman classes
(see \cite{BM2}). Quasianalytic functions are smooth, and, by
Lemma \ref{L:dc}, among the quasianalytic functions of one real
variable there exist functions which coincide with analytic
functions on at most a discrete set. This property is fundamental
in our constructions.
\par In \cite{CLP} we proved that quasianalytic {\it curves} are
pluripolar. Raising the dimension of our submanifolds to at least
$n$ we obtain nonpluripolar examples.
\par Using the same ideas as in the proof of Theorem \ref{T:tor},
we construct in Theorem \ref{T:real} a quasianalytic function on
${\Bbb R}^n$ whose graph intersects any real analytic curve in a
discrete set and, consequently, does not contain any analytic
curve. It serves as an example of an ``extremely'' smooth function
which is not {\it arc-analytic} anywhere (see \cite{BM1}).
\par We would like to thank Al Taylor for introducing us to
the notion of quasianalytic functions.

\section{Basic definitions and facts}\label{S:bdf}
\par Let $f:[a,b]\ra{\Bbb R}$ be a $C^\infty$ function and let
$$M_j(f)=\sup_{a\le x\le b}|f^{(j)}(x)|.$$
\par Given an increasing sequence $\{M_j\}$ which is logarithmically
convex, the class $C^\#\{M_j\}$ consists of all smooth functions
$f:[a,b]\ra{\Bbb R}$ satisfying the estimate $M_j(f)\leq C^jM_j$
for all $j$, with a constant $C$ depending on $f$. Note that if
$f,\,g\in C^\#\{M_j\}$ then clearly $f+g\in C^\#\{M_j\}$.
\par The proof of the following lemma can be found in
\cite[Prop. 8.4.1]{H}.
\bL\label{L:caf} Let $\{M_j\}$ be a sequence such that $M_j^{1/j}>j$
is increasing and let $f\in C^\#\{M_j\}$ on some interval $[a,b]$.
If $g:[\alpha,\beta]\ra[a,b]$ is a real analytic function, then
the function $h=f\circ g\in C^\#\{M_j\}$ on $[\alpha,\beta]$.\eL

\par Let $$\tau(r)=\inf_{j\geq1}\frac{M_j}{r^j}$$ be the associated
function to the sequence $\{M_j\}$. We denote by $\nu(r)$ the
largest integer such that $\tau(r)=M_{\nu(r)}/r^{\nu(r)}$. The
following lemma is a slight extension of a result by S.
Mandelbrojt in \cite[Ch.VI.41]{M}.

\bL\label{L:dc} Let $\{M_j\}$ be a sequence such that $C^j=o(M_j)$
for every constant $C>1$. Suppose
$$\liminf_{r\to\infty}\frac{\nu(r)}{\log r}\ge A>0.$$
There exists a smooth periodic function $f\in C^\#\{M_j\}$ on
$[0,2\pi]$ with the following property: for every sequence
$\{N_j\}$ such that
$$\lim_{j\to\infty}\frac{{N_j}^{1/j}}{M_j^{1/j}}=0,$$
$f$ does not belong to $C^\#\{N_j\}$ on any interval.\eL

\begin{pf} Following \cite[Ch.VI.41]{M}, we introduce
$$a_n=\frac{\tau(2^n)}{2^n}$$ and for $k=1,2,\dots$ we define
$$f_k(x)=\sum_{n=k}^\infty a_n\cos 2^nx.$$ Since $\tau(2^n)\le
M_j2^{-jn}$ for all $j=1,2,\dots$, we have $a_n\le M_j2^{-(j+1)n}$
and
$$|f^{(j)}_k(x)|\le\sum_{n=k}^\infty a_n2^{jn}\le
M_j\sum_{n=k}^\infty 2^{-n}\le M_j.$$ Thus each of the functions
$f_k\in C^\#\{M_j\}$ on $[0,2\pi]$.
\par We show that $f=f_1$ is the desired function.
Suppose that $f_1$ is in $C^\#\{N_j\}$ on some interval
$[\al,\beta]$ in $[0,2\pi]$. We may assume that $\al=2\pi
m2^{-k_0}$ and $\beta=2\pi (m+1)2^{-k_0}$ for some appropriate
choice of integers $m$ and $k_0$ and we write
$f_1=g_{k_0}+f_{k_0}$, where
$$g_{k_0}(x)=\sum_{n=1}^{k_0-1}a_n\cos 2^nx.$$ We have
$|g^{(j)}_{k_0}(x)|\le C^j$ for all $j=1,2,\dots$ and
$x\in[0,2\pi]$, where $C$ is a constant depending on $k_0$. Since
$$\lim_{j\to\infty}\frac{N_j^{1/j}}{M_j^{1/j}}=0,$$
we conclude that the function $f_{k_0}\in C^\#\{L_j\}$ on
$[\al,\beta]$, where $L_j=N_j+C^j$ and
$$\lim_{j\to\infty}\frac{L_j^{1/j}}{M_j^{1/j}}=0.$$ However, this
function is periodic with period $2^{-k_0+1}\pi$. Hence
$f_{k_0}\in C^\#\{L_j\}$ on $[0,2\pi]$.
\par Since
$$f^{(j)}_{k_0}(x)=\sum_{n=k_0}^\infty a_n2^{jn}\phi_n(x)$$ where
$\phi_n(x)$ is either $\pm\cos 2^nx$ or $\pm\sin 2^nx$,
$$|a_n|=
\left|\frac1{\pi2^{nj}}\il0^{2\pi}f^{(j)}_{k_0}(x)\phi_n(x)\,dx\right|
\le\frac{2R^jL_j}{2^{nj}}\;.$$ Here $R\geq 1$ depends on $k_0$.
\par If $j_n=\nu(2^n)$, then $a_n=M_{j_n}2^{-(j_n+1)n}$ by the definition
of $a_n$. By the previous inequality, this implies  that
$$M_{j_n}^{1/j_n}\le C_1L_{j_n}^{1/j_n}2^{n/j_n}$$
or
$$\lim_{n\to\infty}\frac{L_{j_n}^{1/j_n}}{M_{j_n}^{1/j_n}}\ge
\frac1{C_1}\liminf_{n\to\infty}2^{-n/j_n}.$$ By hypothesis
$$\liminf_{n\to\infty}\frac{j_n}{n}\ge B=A\log2>0;$$ thus we get
$$\lim_{n\to\infty}\frac{L_{j_n}^{1/j_n}}{M_{j_n}^{1/j_n}}\ge
\frac{2^{-1/B}}{C_1}>0.$$ This contradiction proves the lemma.
\end{pf}

\vskip 2mm The class $C^\#\{M_j\}$ on $[a,b]$ is called
{\it quasianalytic} if
every function in $C^\#\{M_j\}$ which vanishes to infinite order
at some point in $[a,b]$ must be identically equal to 0. A smooth
function $f$ is called {\it quasianalytic in the sense of Denjoy}
if the class $C^\#\{M_j(f)\}$ is quasianalytic. The
Denjoy-Carleman theorem (see e.g. \cite[3.10(12)]{Ti}) states that
the class $C^\#\{M_j\}$ is quasianalytic if
$$\sum_{j=1}^\infty\frac1{M_j^{1/j}}=\infty.$$
\section{Examples}
\par We introduce the notation
$$\log_1x=\log x,\;\log_kx=\log_{k-1}(\log x),$$ and
consider the sequences $\{M^{(k)}_j\}$, where
$M^{(k)}_j=(j\log_kj)^j$ and $k>0$. These sequences are
logarithmically convex and it is easy to see that they verify the
hypotheses of Lemmas \ref{L:caf} and \ref{L:dc}. Moreover, the
classes $C^\#\{M^{(k)}_j\}$ are quasianalytic.
\par Our constructions rely on the following examples of smooth
hypersurfaces in ${\Bbb R}^m$, diffeomorphic to ${\Bbb R}^{m-1}$,
which intersect any real analytic curve in a discrete set. Let
$f_k:{\Bbb R}\ra{\Bbb R}$, $k=1,\dots,m-1$, be smooth periodic
functions such that $f_k\in C^\#\{M^{(k)}_j\}\sm
C^\#\{M^{(k+1)}_j\}$ on any interval. Such functions exist by
Lemma \ref{L:dc}.

\bT\label{T:real} The hypersurface $H\subset{\Bbb R}^m$, defined
by the equation $$x_m=\sum_{j=1}^{m-1}f_j(x_j),$$ intersects any
real analytic curve in a discrete set. \eT

\begin{pf} Let $\gamma(t)=(\gamma_1(t),\dots,\gamma_m(t))$ be a
real analytic curve defined for $t$ in an open interval $I$ about
$0$ with $\gamma(0)=x_0=(x_{01},\dots,x_{0m})\in H$ and
such that there exists a sequence $t_s\ra0$, $t_s\neq0$, with
$\gamma(t_s)\in H$. It suffices to show that $\gamma(t)=x_0$ for
all $t\in I$.
\par We have
\begin{equation}\label{e:id}
\gamma_m(t_s)=\sum_{j=1}^{m-1}f_j(\gamma_j(t_s)),\;\forall\,s\geq1.
\end{equation}
By Lemma \ref{L:caf}, the right side is quasianalytic; since the
left side is real analytic, it follows that (\ref{e:id}) holds for all
$t\in I$. We show by induction on $k=1,\dots,m-1$ that all
functions $\gamma_k$ are identically constant. Then $\gamma_m$ is
also identically constant by (\ref{e:id}), and we are done.
\par If the function $\gamma_1$ is not
identically constant, then it has a real analytic inverse function on
some interval $J$. For $u\in J$, we have by (\ref{e:id})
$$f_1(u)=\gamma_m(\gamma_1^{-1}(u))-
\sum_{j=2}^{m-1}f_j(\gamma_j(\gamma_1^{-1}(u))).$$ By Lemma
\ref{L:caf} and the choice of $f_1\in C^\#\{M^{(1)}_j\}\sm
C^\#\{M^{(2)}_j\}$ this is impossible. Hence $\gamma_1(t)\equiv
x_{01}$. We assume by induction that for $2\leq k\leq m-1$ we have
$\gamma_j(t)\equiv x_{0j}$ for all $j\leq k-1$. Equation
(\ref{e:id}) becomes
$$\gamma_m(t)-\sum_{j=1}^{k-1}f_j(x_{0j})=
\sum_{j=k}^{m-1}f_j(\gamma_j(t)).$$ A similar argument shows,
since $f_k\in C^\#\{M^{(k)}_j\}\sm C^\#\{M^{(k+1)}_j\}$, that
$\gamma_k(t)\equiv x_{0k}$.\end{pf}

\bT\label{T:me} There exists a smooth submanifold $R$ in ${\Bbb
C}^n$, diffeomorphic to ${\Bbb R}^{2n-2}$, which intersects every
analytic disk in a discrete set. Moreover, $R$ contains a smooth
submanifold $M$, diffeomorphic to ${\Bbb R}^n$, which is not
pluripolar.\eT
\begin{pf}
Using Lemma \ref{L:dc}, we choose smooth periodic functions
$f_k,\,g_k:{\Bbb R}\ra{\Bbb R}$, $k=1,\dots,n-1$, such that
$$f_k\in C^\#\{M^{(2k-1)}_j\}\sm
C^\#\{M^{(2k)}_j\},\;g_k\in C^\#\{M^{(2k)}_j\}\sm
C^\#\{M^{(2k+1)}_j\},$$ on any interval.
\par Suppose that in ${\Bbb C}^n\equiv{\Bbb R}^{2n}$ we have coordinates
$$z=(z_1,z_2,\dots,z_n),\;z_j=x_j+iy_j=(x_j,y_j).$$
Consider the real hypersurface $$S=\{z\in{\Bbb
C}^n:\,y_n=F(z_1,\dots,z_{n-1},x_n)\},$$ where $F:{\Bbb
R}^{2n-1}\ra{\Bbb R}$ is a stricly convex real analytic function
(i.e., its Hessian is positive definite at any point). We define
the submanifold $R\subset S$ by the equation
$$x_n=\sum_{j=1}^{n-1}(f_j(x_j)+g_j(y_j)).$$

\par Clearly $R$ is diffeomorphic to ${\Bbb R}^{2n-2}$. Suppose,
for the sake of obtaining a contradiction, that $z_0\in R$,
$\phi=(\phi_1,\dots,\phi_n):\,U\to{\Bbb C}^n$ is an analytic disk
with $\phi(0)=z_0$, and there is a sequence $\zeta_s\to0$ in $U$
such that the points $\phi(\zeta_s)\in R$ are distinct. The set of
$\zeta=x+iy\in U$ with $\phi(\zeta)\in S$ is defined by the
equation $\Phi(x,y)=0$, where
$$\Phi(x,y)=\im\phi_n(\zeta)-
F(\phi_1(\zeta),\dots,\phi_{n-1}(\zeta),\re\phi_n(\zeta)).$$ The
function $\Phi$ is real analytic and not identically 0, since $S$
does not contain any analytic disk.
\par It follows from the Weierstrass Preparation Theorem and Theorem
3.2.5 in \cite{KP} that there are integers $k,\,l>0$ such that the
solutions of either of the equations
$\Phi(t^k,y)=0,\,\Phi(-t^l,y)=0$ in some neighborhood of the
origin are graphs of finitely many real analytic functions
$y=\al(t)$, defined in open intervals about $0$, or $\{t=0\}$ (see
also \cite{BM1}). We conclude that there exist a sequence
$t_s>0,\,t_s\ra0$, and a real analytic curve
$\gamma(t)=(\gamma_1(t),\gamma_2(t))$ defined for $t$ in an
interval about 0 (where $\gamma_1(t)$ is one of the functions
$t^k$, $-t^l$, or $0$), such that $\gamma(0)=0$ and
$$\phi(\gamma_1(t_s)+i\gamma_2(t_s))\in R,\;\forall\,s\geq1.$$
By passing to a subsequence, we can assume
$\zeta_s=\gamma_1(t_s)+i\gamma_2(t_s)$.
\par If $\psi_j(t)=\re\phi_j(\gamma_1(t)+i\gamma_2(t)),\;
\eta_j(t)=\im\phi_j(\gamma_1(t)+i\gamma_2(t))$, for $j=1,\dots,n$,
we have $$\psi_n(t_s)=
\sum_{j=1}^{n-1}\left(f_j(\psi_j(t_s))+g_j(\eta_j(t_s))\right).$$
Hence the hypersurface $H\subset{\Bbb R}^{2n-1}$ defined by the
equation $$x_n=\sum_{j=1}^{n-1}(f_j(x_j)+g_j(y_j))$$ intersects
the real analytic curve $$\Psi(t)=
(\psi_1(t),\eta_1(t),\dots,\psi_{n-1}(t),\eta_{n-1}(t),\psi_n(t))$$
at the points $\Psi(t_s)$. By Theorem \ref{T:real} the functions
$\psi_n,\,\psi_j,\,\eta_j$, $1\leq j\leq n-1$, are constant, so
$\phi(\zeta_s)=(\eta,\beta +i\delta_s)$ for some $\eta\in{\Bbb
C}^{n-1}$, $\beta\in{\Bbb R}$, and a real sequence $\delta_s$.
Since $\phi(\zeta_s)\in S$ we have $\phi(\zeta_s)=z_0$, a
contradiction. This shows that $R$ intersects any analytic disk in
a discrete set.

We conclude the proof with the construction of a submanifold $M$
of $R$ with the desired properties. We choose $F$ so that it has a
minimum point at the origin and $F(0)=0$. The functions
$f_j,\,g_j$ we select from Lemma \ref{L:dc} can be chosen so that
$f_j(0)=g_j(0)=0$ and $f_1'(0)=1$, $g_1'(0)=0$. Let $M\subset R$
be defined by $y_2=\dots=y_{n-1}=0$. Clearly $M$ is diffeomorphic
to ${\Bbb R}^n$ and $0\in M$. The tangent space $T_0M$ is given by
the equations $y_2=\dots=y_n=0$,
$x_n=x_1+\sum_{j=2}^{n-1}f_j'(0)x_j$, so $T_0M\cap iT_0M=\{0\}$.
Hence $M$ is not pluripolar, since it is generating at $0$ (see
\cite{P} and \cite{Sa}).
\end{pf}

\vskip 2mm We now proceed with the proof of Theorem \ref{T:tor}
stated in the Introduction. We need the following lemma.

\bL\label{L:inv} Let $G$ and $H$ be real valued, real analytic functions,
defined in open intervals about $\theta_0$ and $0$.
Assume that $k<\infty$ is the vanishing order of $G-G(\theta_0)$
at $\theta_0$ and that $H(t_j)=G(\theta_j)$, where
$t_j>0,\,t_j\ra0$, $\theta_j\ra\theta_0$. Then there is an
analytic function $h$ defined in an open interval about $0$ and a
subsequence $\{j_n\}$ so that $\theta_{j_n}=h(t_{j_n}^{1/k})$.\eL

\begin{pf} Without loss of generality we can assume
$\th_0=G(\th_0)=0$. We write $G(\th)=\pm\th^kG_1(\th)$, where
$G_1(0)>0$. Since $|\th_jG_1^{1/k}(\th_j)|=|H(t_j)|^{1/k}$, there
exists a choice of signs such that for infinitely many $j$
$$\th_jG_1^{1/k}(\th_j)=\pm(\pm H(t_j))^{1/k},$$ where $\pm
H(t_j)>0$. Let $H_1(t)$ be the analytic branch of the function
$(\pm H(t^k))^{1/k}$ which is positive for $t>0$. Then
$\th_jG_1^{1/k}(\th_j)=\pm H_1(t_j^{1/k})$. Finally, let
$h(t)=g(\pm H_1(t))$, where $g$ is the inverse of the analytic
function $\th\ra\th G_1^{1/k}(\th)$.
\end{pf}
\noindent{\it Proof of Theorem \ref{T:tor}.} Let $\{a_k\}$ be any
sequence such that $a_1>0$ and $a_1+\dots+a_k<a_{k+1}$ holds for
all $k\geq1$. We define inductively maps
$$T^k:{\Bbb R}^k\ra{\Bbb R}^{k+1},\;T^k=(T^k_1,\dots,T^k_{k+1}),$$
whose image is a $k$-dimensional torus embedded in ${\Bbb
R}^{k+1}$. Let
$$T^1(\theta_1)=(T^1_1(\theta_1),T^1_2(\theta_1))=
(a_1\sin\theta_1,a_1\cos\theta_1).$$ Given $T^k$ we define
$T^{k+1}$ by
\begin{eqnarray*}
T^{k+1}_j(\theta_1,\dots,\theta_k,\theta_{k+1}) & = &
T^k_j(\theta_1,\dots,\theta_k),\;\;\;\;\;\;\;j=1,\dots,k,\\
T^{k+1}_{k+1}(\theta_1,\dots,\theta_k,\theta_{k+1}) & = &
(T^k_{k+1}(\theta_1,\dots,\theta_k)+a_{k+1})\sin\theta_{k+1},\\
T^{k+1}_{k+2}(\theta_1,\dots,\theta_k,\theta_{k+1}) & = &
(T^k_{k+1}(\theta_1,\dots,\theta_k)+a_{k+1})\cos\theta_{k+1}.
\end{eqnarray*}
It follows by induction on $k$ that
$|T^k_j(\theta_1,\dots,\theta_k)|\leq a_1+\dots+a_k$ for $j\leq
k+1$, and that $T^k$ is injective on $[0,2\pi)^k$.

We use notation similar to that used in the proof of Theorem \ref{T:me}. Let
$f_k:{\Bbb R}\ra{\Bbb R}$, $k=1,\dots,2n-2$, be smooth periodic
functions such that $f_k\in C^\#\{M^{(k)}_j\}\setminus
C^\#\{M^{(k+1)}_j\}$ on any interval. Let $R$ be the submanifold
of a sphere $S\subset{\Bbb C}^n$ of sufficiently large radius $r$,
given by the image of the map $z=F(\theta_1,\dots,\theta_{2n-2})$
defined by
\begin{eqnarray}
 & & x_j=T^{2n-2}_{2j-1}(\theta_1,\dots,\theta_{2n-2}),\;
 y_j=T^{2n-2}_{2j}(\theta_1,\dots,\theta_{2n-2}),\;j\leq n-1,
 \label{e:j}\\
 & & x_n=T^{2n-2}_{2n-1}(\theta_1,\dots,\theta_{2n-2})+
 \sum_{j=1}^{2n-2}f_j(\theta_j),\label{e:xn}\\
 & & y_n=(r^2-x_1^2-y_1^2-\dots-x_{n-1}^2-y_{n-1}^2-x_n^2)^{1/2}.
 \label{e:yn}
\end{eqnarray}

Let $X$ denote the image of the map
$z=G(\theta_1,\dots,\theta_{2n-2})$ defined by
\begin{eqnarray}
x_j & = & \left(\sum_{k=1}^{2j-2}a_k\cos\theta_k\dots
\cos\theta_{2j-2}+a_{2j-1}\right)\sin\theta_{2j-1},\label{e:txj}\\
y_j & = & \left(\sum_{k=1}^{2j-1}a_k\cos\theta_k\dots
\cos\theta_{2j-1}+a_{2j}\right)\sin\theta_{2j},\label{e:tyj}\\
x_n & = & \left(\sum_{k=1}^{2n-3}a_k\cos\theta_k\dots
\cos\theta_{2n-3}+a_{2n-2}\right)\cos\theta_{2n-2},\label{e:txn}\\
y_n & = &
\left(r^2-x_1^2-y_1^2-\dots-x_{n-1}^2-y_{n-1}^2-x_n^2\right)^{1/2}.
\label{e:tyn}
\end{eqnarray}
Here $1\leq j\leq n-1$ and we used the explicit formulas defining
$T^{2n-2}$. Then $X$ is diffeomorphic to the $(2n-2)$-dimensional
torus. By choosing $f_j$ with sufficiently small $C^1$-norm, the
map $F$ is a $C^1$-perturbation of the map $G$, hence $R$ is also
diffeomorphic to the $(2n-2)$-dimensional torus.

\par Suppose, for the sake of obtaining a contradiction, that there are $z_0\in R$ and a
non-constant analytic disk $\phi=(\phi_1,\dots,\phi_n):U\ra{\Bbb
C}^n$ with $\phi(0)=z_0$ and with the following property: there
exist sequences of points $\zeta_s\ne0$, $\zeta_s\ra0$ in $U$, and
$\Th_s=(\th_{1s},\dots,\th_{(2n-2)s})\ra\Th_0=
(\th_{10},\dots,\th_{(2n-2)0})$ in ${\Bbb R}^{2n-2}$, such that
$F(\Th_s)=\phi(\zeta_s)$.

\par The set $\{\zeta=x+iy\in U:\,\phi(\zeta)\in S\}$ is defined
by the equation
$$\Phi(x,y):=\sum_{j=1}^n|\phi_j(x+iy)|^2-r^2=0.$$
The function $\Phi$ is real analytic and is not identically equal
to 0, since $S$ does not contain any analytic disk. Hence there
exist a real analytic curve $\gm$ defined in an interval about 0, with
$\gm(0)=0$, and a sequence $t_s>0,\,t_s\ra0$, such that
$\Phi(\gm(t))\equiv0$ and $\gm(t_s)=\zeta_s$ (after passing to a
subsequence of $\{\zeta_s\}$, if necessary).

We let $\psi_j(t)=\re\phi_j(\gamma(t))$,
$\eta_j(t)=\im\phi_j(\gamma(t))$, $j=1,\dots,n$, and we write the
equations (\ref{e:j})-(\ref{e:xn}) corresponding to
$\phi(\gamma(t_s))=F(\Theta_s)$. Equation (\ref{e:j}) for $x_1$
becomes $\psi_1(t_s)=a_1\sin\theta_{1s}$. By Lemma \ref{L:inv},
there exist an integer $N_1\geq1$ and a real analytic function $h_1$
near 0 so that $\theta_{1s}=h_1(t_s^{1/N_1})$ for a sequence of
integers $s\ra\infty$. Changing variables $t=u_1^{N_1}$ and
letting $u_{1s}=t_s^{1/N_1}$, we obtain from equation (\ref{e:j})
for $y_1$ that
$$\eta_1(u_{1s}^{N_1})=(a_1\cos(h_1(u_{1s}))+a_2)\sin\theta_{2s}.$$
Now Lemma \ref{L:inv} implies the existence of an integer
$N_2\geq1$ and a real analytic function $h_2$ near 0 such that
$\theta_{2s}=h_2(u_{1s}^{1/N_2})$ for infinitely many $s$. Next we
change variables $t=u_1^{N_1}=u_2^{N_1N_2}$,
$u_{2s}=t_s^{1/(N_1N_2)}$ and consider equation (\ref{e:j}) for
$x_2$:
$$\psi_2(u_{2s}^{N_1N_2})=(a_1\cos(h_1(u_{2s}^{N_2}))\cos(h_2(u_{2s}))
+a_2\cos(h_2(u_{2s}))+a_3)\sin\theta_{3s}.$$

Continuing like this, we conclude that there are an integer
$N\geq1$, real analytic functions $g_j$ near 0, $j=1,\dots,2n-2$, and a
sequence of integers $s\ra\infty$, such that
$\theta_{js}=g_j(u_s)$, where $u_s=t_s^{1/N}$. Using (\ref{e:xn})
we get
$$\psi_n(u_s^N)-T^{2n-2}_{2n-1}(g_1(u_s),\dots,g_{2n-2}(u_s))=
\sum_{j=1}^{2n-2}f_j(g_j(u_s)).$$ Hence the hypersurface
$H\subset{\Bbb R}^{2n-1}$ defined by the equation
$$\theta_{2n-1}=\sum_{j=1}^{2n-2}f_j(\theta_j)$$ intersects the
real analytic curve
$$\Psi(u)=\left(g_1(u),\dots,g_{2n-2}(u),\psi_n(u^N)-
T^{2n-2}_{2n-1}(g_1(u),\dots,g_{2n-2}(u))\right)$$ at the points
$\Psi(u_s)$. By Theorem \ref{T:real} all the functions $g_j$ are
identically constant; thus $F(\Theta_s)=\phi(\zeta_s)=z_0$, a
contradiction. Hence $R$ intersects any analytic disk in a
discrete set.

\par We proceed now with the construction of a submanifold $M$ of
$R$ with the desired properties. Recall that $R$ was defined as
the image of the map $F$ in (\ref{e:j})-(\ref{e:yn}), which is a
$C^1$-perturbation of the map $G$ given by
(\ref{e:txj})-(\ref{e:tyn}). Let $M\subset R$ and $Y\subset X$ be
the submanifolds defined by $x_1=\dots=x_{n-2}=0$, i.e., by taking
$\theta_{2j-1}=0$, $j=1,\dots,n-2$, in the formulas
(\ref{e:j})-(\ref{e:yn}) and in (\ref{e:txj})-(\ref{e:tyn}). Since
$Y$ is diffeomorphic to the $n$-torus, so is $M$. We check that
$Y$ is generating at the point $P$ corresponding to
$\theta_{2k}=0$, $k=1,\dots,n-2$, $\theta_{2n-3}=0$ and
$\theta_{2n-2}=\pi/2$. Indeed, $P$ has ${\Bbb R}^{2n}$-coordinates
$P=(0,\dots,0,a,0,b)$ for some $a,b$, and the tangent space $T_PY$
is given by $x_1=\dots=x_{n-2}=y_{n-1}=y_n=0$. We conclude that
$M$ is generating at the point $P'$ corresponding to the same
values of parameters, so that $M$ is not pluripolar (see \cite{P},
\cite{Sa}). $\Box$

\vskip 2mm\noindent {\bf Remark.} Let $D\subset{\Bbb C}^n$ be any
bounded pseudoconvex domain with real analytic boundary. By
\cite{DF} the boundary of $D$ does not contain any non-constant
analytic disk. The construction in the proof of Theorem
\ref{T:tor} shows the existence of smooth nonpluripolar compact
submanifolds $M$ and $R$ of $\partial D$ which are diffeomorphic
to an $n$-dimensional torus and a $(2n-2)$-dimensional torus, and
which intersect any analytic disk in a discrete set.


\begin{thebibliography}{99999}
\bibitem[B]{B} E. Bedford, {\em Survey of pluripotential theory}, in
{\em Several Complex Variables: Proceedings of the Mittag-Leffler
Institut 1987-1988}, J.-E. Forn\ae ss (ed.) Math. Notes {\bf 38}
Princeton Univ. Press, 1993, 48--97.
\bibitem[BM1]{BM1} E. Bierstone and P. Milman, {\em Arc-analytic
functions}, Invent. Math. {\bf 101} (1990), 411--424.
\bibitem[BM2]{BM2} E. Bierstone and P. Milman, {\em Resolution of
singularities in Denjoy--Carleman classes,} Selecta Math., (to
appear).
\bibitem[CLP]{CLP} D. Coman, N. Levenberg and E. A. Poletsky, {\em
Quasianalyticity and pluripolarity}, (preprint).
\bibitem[DF]{DF} K. Diederich and J. E. Forn\ae ss, {\em Proper
holomorphic maps onto pseudoconvex domains with real-analytic
boundary}, Annals of Math. {\bf 110} (1979), 575--592.
\bibitem[H]{H} L. H\"ormander, {\em The Analysis of Linear Partial Differential Operators I,}
Springer-Verlag, 1990.
\bibitem[Ki]{Ki} C. O. Kiselman, {\em The growth of restriction of
plurisubharmonic functions, Math. Analysis and Applications}, Part
B. Adv. Math., Suppl. Studies 7B (1981), 435--454.
\bibitem[KP]{KP} S. G. Krantz and H. R. Parks, {\em A Primer of
Real Analytic Functions}, Birkh\"auser, 1992.
\bibitem[M]{M} S. Mandelbrojt, {\em S\'eries de Fourier et classes
quasi-analytique de fonctions}, Paris, Gauthier-Villars, 1935.
\bibitem[P]{P} S. I. Pin\v cuk, {\em A boundary
uniqueness theorem for holomorphic functions of several complex
variables}, Mat. Zametki {\bf 15} (1974), 205--212.
\bibitem[R]{R} L. I. Ronkin, {\em Introduction to the theory
of entire functions of several variables}, Amer. Math. Soc.,
Providence, R.I., 1974.
\bibitem [Sa]{Sa} A. Sadullaev, {\em A boundary uniqueness theorem
in ${\Bbb C}^n$,} Mat. Sb. (N.S.) {\bf 101}(143) (1976), 568--583.
\bibitem[T]{Ti} A. F. Timan, {\em Theory of Approximation of Functions
of a Real Variable}, Pergamon Press, Macmillan, New York, 1963.
\end{thebibliography}
\end{document}